\documentclass[aof]{imsart}
\usepackage[dvipsnames]{xcolor}
\usepackage{pdfpages}
\usepackage{amsmath}
\usepackage{amsfonts}
\usepackage{upgreek}
\usepackage{sgame}
\usepackage{accents}
\usepackage{amssymb}
\usepackage{tikz-cd}
\usepackage{float}
\usepackage[numbers]{natbib}   
\usepackage{dutchcal}
\usepackage{pdfpages}
\usepackage{amsmath}
\usepackage{amsfonts}
\usepackage{upgreek}
\usepackage{sgame}
\usepackage{mathtools}
\usepackage{amsthm}
\usepackage{amssymb}
\usepackage{enumitem}
\usepackage{epigraph}
\usepackage{sectsty}
\usetikzlibrary{calc}

\newtheorem{theorem}{Theorem}[section]
\newtheorem{proposition}[theorem]{Proposition}
\newtheorem{lemma}[theorem]{Lemma}
\newtheorem{corollary}[theorem]{Corollary}
\newtheorem{claim}[theorem]{Claim}

\theoremstyle{definition}

\newtheorem{definition}[theorem]{Definition}
\newtheorem{example}[theorem]{Example}

\newtheorem{question}[theorem]{Question}

\newtheorem{remark}[theorem]{Remark}

\newtheorem{aside}[theorem]{Aside}
\newtheorem{note}[theorem]{Note}

\usepackage{hyperref}
\hypersetup{
    colorlinks=true,
    linkcolor=Mulberry,
    filecolor=Mulberry,      
    urlcolor=Mulberry,
    citecolor=Mulberry,
}


\begin{document}
\title{A Simple Proof of the Monotonicity of the Invariant Distribution for a Discrete Markov Chain}
\author{Mark Whitmeyer\thanks{Department of Economics, University of Texas at Austin. Email: \href{mailto:mark.whitmeyer@utexas.edu}{mark.whitmeyer@utexas.edu}}}
\date{\today{}}
\begin{abstract}
This note presents a simple proof of the monotonicity of the invariant distribution of a discrete Markov chain with a finite state space. This answers a question recently raised by David Siegmund.
\end{abstract}


\maketitle

\section{Introduction}

In the 2018 Symposium on Optimal Stopping at Rice University (in memory of Larry Shepp), David Siegmund asked whether there is a simple proof of the following result.

\begin{theorem}
Let $M = (S,P)$ be a (Markov) model with a finite state space $S$, transition matrix $P = \big\{p(i,j)\big\}$ and limit (invariant) distribution $\pi > 0$. Let $M'(S,P')$ be a model with invariant distribution $\pi'$ with (perturbed) matrix $P' = \big\{p'(i,j)\big\}$ such that for some state $s_{0} \in S$,

\[p'(i,s_{0}) \geq p(i,s_{0}), \hspace{1cm} p'(i,j) \leq p(i,j)\]
for all $i$ and for all $j\neq s_{0}$ and at least one $p'(i,s_{0}) > p(i,s_{0})$. Then $\pi'(s_{0}) > \pi(s_{0})$.
\end{theorem}

In recent work, Isaac M. Sonin provides an alternative proof of this using the idea of a censored Markov chain \cite{son}. Here, these techniques are not used, and instead the result is obtained through properties of the expected first return time.
\begin{proof}

Let the number of states be $n$. Without loss of generality set $s_{0} = 1$. We state the following standard results. Let $\mu_{1}$ be the expected first return time to state $1$. Then,

\[\tag{$1$}\label{1}\pi_{1} = \frac{1}{\mu_{1}}\]
Let $\mu_{ij}$ be the expected first hitting time to $j$ from state $i$. Then,
\[\tag{$2$}\label{2}\mu_{1} = 1 + \sum_{k=2}^{n}p(1,k)\mu_{k1}\]
\[\tag{$3$}\label{3}\mu_{j1} = 1 + \sum_{k=2}^{n}p(j,k)\mu_{k1} \hspace{1cm} \text{for} \hspace{2mm} j \geq 2\]
\begin{remark}
It is sufficient to show that for a matrix $P$, the invariant probability $\hat{\pi}_{1} > \pi_{1}$ for a perturbed matrix $\Hat{P} = \big\{\Hat{p}(i,j)\big\}$ where \[\Hat{p}(i,1) = p(i,1) + c_{i}\]
and \[\Hat{p}(i,2) = p(i,2) - c_{i}\]
for all $i$ and for all (feasible) $c_{i} \geq 0$ with at least one $c_{i} > 0$ and $\Hat{p}(i,j) = p(i,j)$ for all $j \neq 1,2$.
\end{remark}
We have for all $i$

\[\frac{\partial{p(i,2)}}{\partial{p(i,1})} = \lim_{c_{i} \to 0} \frac{p(i,2) - c_{i} - p(i,2)}{c_{i}} = -1\] Then, from Expression \ref{2}, we have 

\[\tag{$4$}\label{4}\frac{\partial{\mu_{1}}}{\partial{p(1,1)}} = -\mu_{2,1} + \sum_{k=2}^{n}p(1,k)\frac{\partial{\mu_{k1}}}{\partial{p(1,1)}}\]
and 
\[\tag{$5$}\label{5}\frac{\partial{\mu_{1}}}{\partial{p(j,1)}} = \sum_{k=2}^{n}p(1,k)\frac{\partial{\mu_{k1}}}{\partial{p(j,1)}}\]
for $j \geq 2$. From Expression \ref{3}, we have
\[\tag{$6$}\label{6}\frac{\partial{\mu_{j1}}}{\partial{p(1,1)}} = \sum_{k=2}^{n}p(j,k)\frac{\partial{\mu_{k1}}}{\partial{p(1,1)}}\]
and 
\[\tag{$7$}\label{7}\frac{\partial{\mu_{j1}}}{\partial{p(j,1)}} = -\mu_{j1} + \sum_{k=2}^{n}p(j,k)\frac{\partial{\mu_{k1}}}{\partial{p(j,1)}}\]
for $j \geq 2$. We may combine Expressions \ref{4} and \ref{6} and iterating forward, obtain

\[\tag{$8$}\label{8}\begin{split}
    \frac{\partial{\mu_{1}}}{\partial{p(1,1)}} &= -\mu_{2,1} + \sum_{k_{1}=2}^{n}p(1,k_{1})\sum_{k_{2}=2}^{n}p(k_{1},k_{2})\cdots = -\mu_{2,1} + 0 < 0
\end{split}\]
In a similar fashion, we combine Expressions \ref{5} and \ref{7} and iterate forward,

\[\tag{$9$}\label{9}\begin{split}
    \frac{\partial{\mu_{1}}}{\partial{p(j,1)}} &= \sum_{k=2}^{n}p(1,k)\bigg(-\mu_{j1} + \sum_{k_{1}=2}^{n}p(k,k_{1})\bigg(-\mu_{j1} + \sum_{k_{2}=2}^{n}p(k_{1},k_{2})\cdots \bigg)\cdots \bigg) < 0
\end{split}\]

Since $\mu_{1}$ is strictly decreasing in $p(i,1)$, $\pi_{1}$ must be strictly increasing in $p(i,1)$, and the result is shown.

\end{proof}

\bibliography{sample.bib}

\begin{thebibliography}{1}
\providecommand{\natexlab}[1]{#1}
\providecommand{\url}[1]{\texttt{#1}}
\expandafter\ifx\csname urlstyle\endcsname\relax
  \providecommand{\doi}[1]{doi: #1}\else
  \providecommand{\doi}{doi: \begingroup \urlstyle{rm}\Url}\fi

\bibitem[Sonin(2018)]{son}
Isaac~M. Sonin.
\newblock The answer for a question of david siegmund (siegmund's
  monotonicity).
\newblock \emph{Mimeo}, 2018.

\end{thebibliography}

\end{document}